\documentclass[11pt,a4paper,fleqn]{article}
\usepackage{amsmath}
\usepackage{amssymb}
\usepackage{amscd}
\usepackage[dvips]{graphicx}

\newtheorem{theorem}{Theorem}

\newtheorem{corollary}[theorem]{Corollary}
\newtheorem{definition}[theorem]{Definition}
\newtheorem{example}[theorem]{Example}
\newtheorem{lemma}[theorem]{Lemma}
\newtheorem{proposition}[theorem]{Proposition}
\newtheorem{remark}[theorem]{Remark}
\newenvironment{proof}{{\bf Proof. }}{\hfill$\rule{1ex}{1ex}$\par\medskip}
\begin{document}
\newcommand{\bt}{\begin{theorem}}
\newcommand{\et}{\end{theorem}}
\newcommand{\bd}{\begin{definition}}
\newcommand{\ed}{\end{definition}}
\newcommand{\bs}{\begin{proposition}}
\newcommand{\es}{\end{proposition}}
\newcommand{\bp}{\begin{proof}}
\newcommand{\ep}{\end{proof}}
\newcommand{\be}{\begin{equation}}
\newcommand{\ee}{\end{equation}}
\newcommand{\ul}{\underline}
\newcommand{\br}{\begin{remark}}
\newcommand{\er}{\end{remark}}
\newcommand{\bex}{\begin{example}}
\newcommand{\eex}{\end{example}}
\newcommand{\bc}{\begin{corollary}}
\newcommand{\ec}{\end{corollary}}
\newcommand{\bl}{\begin{lemma}}
\newcommand{\el}{\end{lemma}}

\title{Oriented Steiner quasigroups}
\author{Izabella Stuhl}
\date{}
\maketitle


\begin{abstract}
We introduce the notion of an oriented Steiner quasigroup and
develop elements of a relevant algebraic apparatus. The approach
is based upon (modified) Schreier-type $f$-extensions for quasigroups
(cf. earlier works \cite{S, NSt, NSt2}) achieved through oriented Steiner triple systems.
This is done in a fashion similar to in \cite{SS} where an analogous construction
was established for loops. As a justification of this concept briefly discuss an application of oriented Steiner triple systems
in cryptography using oriented Steiner quasigroups.
\end{abstract}

\section{Introduction}

A {\it Steiner triple system} $\mathfrak{S}$ is an incidence
structure consisting of points and blocks such that every two
distinct points are contained in precisely one block and any block
has precisely three points. A finite Steiner triple system with
$n$ points exists if and only if $n \equiv 1$ or $3$ $(mod\hskip
4pt 6)$.

A possible way to enrich Steiner triple systems is to introduce in each block
a cyclic order.
These objects occur in the literature as oriented Steiner
triple systems $(\mathfrak{S}, T)$ (see \cite{W}) or as Mendelsohn triple systems
(see \cite{M} Section 28, p. 388).

Steiner triple systems play an important role in applications, particulary in cryptography,
coding theory, statistics and computer science. Cf. \cite{CD}, \cite{DS}
(where a number of chapters are dedicated to this concept),
and also \cite{CR} and the references therein. A particular example of an application
is given in \cite{V}. The results of the present paper develop an algebraic background
that can be useful, e.g., in the context of \cite{Be}.

A Steiner triple system $\mathfrak{S}$
provides a multiplication on the set of pairs of different points $x, y$ taking as product the
third point of the block joining $x$ and $y$. Defining $x\cdot x=x$, we get the Steiner quasigroup
associated with $\mathfrak{S}$. Conversely, a Steiner quasigroup determines a Steiner triple
system whose points are the elements of the quasigroup, and
the blocks are the triples $\{x, y, xy\}$ for $x\neq y\in \mathfrak{S}$.

Steiner triple systems and Steiner quasigroups are in a one-to-one correspondence,
thus algebraic structures corresponding to the oriented systems are also of interest.
We also note that a class of idempotent semi-symmetric quasigroups (also called Mendelsohn quasigroups)
can be related to
Mendelsohn triple systems (see e.g. Theorem 28.8 in \cite{M} p. 389.)

In an earlier work \cite{SS} one studied loops corresponding to
oriented Steiner triple systems obtained by a Schreier-type extension
process (developed for loops in \cite{NS2}).
In this paper we determine classes of quasigroups related to oriented
Steiner triple systems in the same way as it was done in \cite{SS}. Here we use
modified Schreier-type $f$-extensions for quasigroups introduced in \cite{S}, \cite{NSt, NSt2}.
Furthermore, we verify which properties are satisfied for these classes of quasigroups, among the properties discussed in \cite{SS} for
the case of loops. Namely, we inspect alternative and flexible laws, inverse and cross
inverse properties, Bol and Moufang identities.
This provides some information about the
similarity and diversity in the behavior of loops and quasigroups corresponding to the
same oriented Steiner triple system.

In the final part of the paper we propose an encryption algorithm based on quasigroup extensions. In our view, this may open
an application of oriented Steiner triple systems
in cryptography, based on properties of oriented Steiner quasigroups.

\section{Preliminaries}

A quasigroup $Q$ satisfies the {\it left inverse property} or the {\it right inverse property}
if there exists a bijection $\iota:Q\to Q$ such that the relation $\iota(x)\cdot(x\cdot y) = y$ or
$(y\cdot x)\cdot\iota(x)= y$ holds for all $x, y\in Q$, respectively.
If both the left and the right inverse properties hold in a quasigroup, then is said to have the {\it inverse property}.
If for any $x\in Q$ there exists an element $x'$ such that $(xy)x'=y$ for all $y\in Q$ then
$Q$ has the {\it left cross inverse property} and if for any $x\in Q$ there exists an element $x''$ such that $x''(yx)=y$ for all $y\in Q$ then
$Q$ has the {\it right cross inverse property}. A quasigroup has the {\it cross inverse property} if satisfies both the left and the right cross inverse properties. $Q$ is {\it left alternative}, respectively {\it right
alternative}
if $ x\cdot(x\cdot y) = x^2\cdot y$, respectively $(y\cdot x)\cdot x = y\cdot x^2$ for all $x, y\in Q$.
A quasigroup $Q$ is {\it flexible} if $ x\cdot(y\cdot x) = (x\cdot y)\cdot x$ for all $x, y\in Q$.
A quasigroup satisfying the identity $x(yx) = y$ is called {\it semi-symmetric}.
$Q$ is a {\it left Bol quasigroup}, respectively {\it right Bol quasigroup} if
it satisfies the identity $x\cdot (y\cdot (x\cdot z)) = (x\cdot(y\cdot x))\cdot z$,
respectively $((x\cdot y)\cdot z)\cdot y = x\cdot ((y\cdot z)\cdot y).$
A quasigroup that satisfies the left and
the right Bol identity is a {\it Moufang quasigroup}. An idempotent totally symmetric quasigroup is
called {\it Steiner quasigroup}. A totally symmetric loop of exponent $2$ is called
{\it Steiner loop}.

An {\it oriented Steiner loop} $L$ is a loop, for which there is an oriented Steiner triple system $(\mathfrak{S},T)$
with the following properties (1)--(3).  (1) $L$ is a loop extension of the group of order $2$ by a Steiner loop $S=\mathfrak{S}\cup e$; (2) the restriction of the factor system of the extension to
$(\mathfrak{S}\times \mathfrak{S})\setminus \{(x,x), x\in
\mathfrak{S}\}$ coincides with the orientation function of
$(\mathfrak{S},T)$; and (3) $f(x,x)=-1$, respectively $f(x,x)=1$ $\forall$
$x\in S\setminus \{e\}$. Cf. Definition 4 in \cite{SS}, p. 136.

Oriented Steiner loops related to oriented Steiner triple systems, in such a way that two non-isomorphic loops are
related to the same oriented Steiner triple system. In a similar manner
we deal with oriented Steiner quasigroups;
to associate a unique quasigroup to an oriented Steiner triple system we also
consider extensions of an object of order $3$.

Let $Q$, $K$ be quasigroups and let $f$
be a function $f:Q\times Q\to K$. Consider the operation
$$
(a,\alpha)\circ (b,\beta):=(ab,f(a,b)\cdot \alpha\beta),
$$
defined on the set $Q\times K=\{(a,\alpha), a\in Q, \alpha \in K\}$.
Since the equations $(a,\alpha)\circ (x,\xi)=(b,\beta)$ and
$(x,\xi)\circ (a,\alpha)=(b,\beta)$ have unique solutions
$\mathcal{Q}_f=(Q\times K, \circ)$ is a quasigroup.
Moreover, the mapping
$(Q\times K,\circ)\to Q:(a,\alpha)\mapsto a$ is a quasigroup homomorphism.
The quasigroup $\mathcal{Q}_f$ is an {\it $f$-extension of
the quasigroup $K$ by the quasigroup $Q$}. The map
$\varphi_{f}:\mathcal{Q}_f\longrightarrow Q:(a,\alpha)\mapsto a$
is the related homomorphism of the $f$-extension
$\mathcal{Q}_f$. The kernel of $\varphi_{f}$ is a quasigroup congruence;
its classes are the equivalence classes of the $f$-extension.

\section{Extensions of commutative quasigroups of order three}\label{ext3}

\bt
Let $(Q,\cdot)$ and $(K,+)$ be two quasigroups. Suppose $(K,\ast)$
is a quasigroup isotopic to $(K,+)$ via a principal
isotopism $(\alpha + \beta)=\varphi_{1}(\alpha)\ast
\varphi_{2}(\beta)$. Let the map $T=(id,\tau): Q\times K\longrightarrow Q\times K$ be defined by
$(a,\alpha)\mapsto (a,\tau(\alpha))$, where $\tau$ is an automorphism of $K$.
\begin{itemize}
\item[(i)] If $\varphi_{1}=\varphi_{2}=\tau$ is an involutory
automorphism, then the quasigroup extensions
$\mathcal{Q}_{f}^{+}:=(Q\times K,+,f)$ and
$\mathcal{Q}_{g}^{\ast}:=(Q\times K,\ast,g)$ are isomorphic if
and only if $f(a,b)=g(a,b)$ for all $a, b\in Q$.
\item[(ii)] For involutory automorphisms
$\varphi_{1}=\varphi_{2}$, $\tau$ the quasigroup extensions
$\mathcal{Q}_{f}^{+}:=(Q\times
K,+,f)$ and $\mathcal{Q}_{g}^{\ast}:=(Q\times K,\ast,g)$ are
isomorphic if and only if $\varphi_{2}(\tau(f(a,b)))=g(a,b)$ for
all $a,b \in Q$.
\end{itemize}
\et

\bp
For the map $T=(id,\tau):\mathcal{Q}_{f}^{+}\longrightarrow \mathcal{Q}_{g}^{\ast}$ we have\\
$$\begin{array}{l}T((a,\alpha)+(b,\beta))=T(ab,(\alpha +\beta)+f(a,b))\\
\quad =(ab,\tau ((\alpha +\beta)+f(a,b)))=
(ab,(\tau(\alpha)+\tau (\beta))+ \tau (f(a,b)))\\
\quad = (ab,
\varphi_{1}(\varphi_{1}(\tau(\alpha))\ast \varphi_{2}(\tau
(\beta)))\ast \varphi_{2}(\tau (f(a,b)))),\end{array}$$ with
$\varphi_{1}=\varphi_{2}$ equal to $(ab, (\tau(\alpha)\ast \tau
(\beta))\ast \varphi_{2}(\tau (f(a,b))))$ and with
$$\varphi_{2}(\tau(f(a,b)))=g(a,b),\;(ab, (\tau(\alpha)\ast \tau
(\beta))\ast g(a,b))= T(a,\alpha)\ast T(b,\beta)$$
for all $a, b\in Q$.

If, in addition $\varphi_{1}=\varphi_{2}=\tau$ then
$g(a,b)=\varphi_{2}(\tau (f(a,b)))=f(a,b)$.
\ep

There are three commutative quasigroups of order $3$, one of which
is the cyclic group $Z_{3}$ of order $3$. The second one is a
commutative quasigroup that has an idempotent element. But then
all three elements are idempotent, and we get the Steiner quasigroup
$K_{3}$ of order $3$. The third quasigroup $Q_{3}$ is a commutative
quasigroup that has no idempotent element. The group $Z_{3}=\Big(\{\alpha,\beta,\gamma=e\},+\Big)$ is
a principal isotope of $K_{3}=\Big(\{\alpha,\beta,\gamma\},\ast\Big)$ and of
$Q_{3}=\Big(\{\alpha,\beta,\gamma\},\diamond\Big)$, with
$\alpha\ast \beta=\varphi_{1}(\alpha)+\varphi_{2}(\beta)=-\alpha -
\beta$ and $\alpha\diamond
\beta=\varphi_{1}(\alpha)+\varphi_{2}(\beta)=-\alpha - \beta -e$, respectively.
\bigskip

Theorem 1 implies
\bigskip

\bc
Quasigroup $f$-extensions of $Z_{3}$ by $Q$, $K_{3}$ by $Q$ and $Q_{3}$ by $Q$
with the same factor system are isomorphic for any quasigroup $Q$.
\ec

\bp
The claim for the first and the second extensions is obtained by choosing equal
involutory automorphisms
$\varphi_{1}=\varphi_{2}=\tau$ given by $\alpha\mapsto
-\alpha $ for all $\alpha$.
In the case of the first and the third extensions the map
$T=(id,\tau):(Q\times Q_{3},\diamond,g)\longrightarrow
(Q\times Z_{3},+,f)$ is an isomorphism with $\tau:\alpha\in
Q_{3} \mapsto -\alpha \in Z_{3}$.
\ep

\section{Quasigroups corresponding to oriented Steiner triple systems}\label{sec4}

With any oriented Steiner triple system $(\mathfrak{S},T)$ one can
associate a function $f^{\ast}:\mathfrak{S}\times \mathfrak{S}\setminus \{(x,x); x\in \mathfrak{S}\}
\longrightarrow \{\pm 1\}$ called the {\it orientation
function} of $(\mathfrak{S},T)$. If $a_{1}$, $a_{2}$ are distinct
points determining the oriented block $(a_{1},a_{2},a_{3})$, then
$f^{\ast}(a_{1},a_{2})=1$ and $f^{\ast}(a_{2},a_{1})=-1$.

Using orientation functions, we determine the factor system
of the $f$-extensions which yield the corresponding quasigroups.

As Steiner quasigroups are commutative inverse property quasigroups,
Steiner quasigroups satisfying left or right Bol identity are
Moufang. According to \cite{K}, quasigroups fulfilling any one of the
Moufang identities have a unit element. In \cite{CLRS}, non-Moufang Steiner loops have been considered, of both kinds: (i) satisfying and (ii) not satisfying Moufang's theorem. Steiner quasigroups do
not possess the alternative law: $y=(yx)x\neq y\cdot x^2=yx$. These facts imply that
quasigroups related to oriented Steiner triple systems via
$f$-extensions by Steiner quasigroups associated with the non oriented
Steiner triple systems cannot be Bol or Moufang quasigroups and cannot have
the alternative law.

Furthermore, the
left, middle and right nuclear square conditions
$$(xx)(yz)=((xx)y)z, \hskip 15pt x((yy)z)=(x(yy))z, \hskip 15pt x(y(zz))=(xy)(zz)$$
in Steiner quasigroups yield the associative law because of idempotency.
But an associative Steiner quasigroup must be trivial, of order $1$. Thus,
the quasigroups obtained as the above $f$-extensions cannot be left, middle or right nuclear square.

\subsection{Oriented Steiner quasigroups}

Similarly to Definition 4 in \cite{SS}, p. 136, we propose
\bd An oriented Steiner quasigroup $\mathcal{Q}_f^{+}$ ( $\mathcal{Q}_f^{-}$) is an $f$-extension
of the cyclic group $Z_2$ of order $2$ by the Steiner quasigroup $Q$ for
which there exists an oriented Steiner triple system
$(\mathfrak{S},T)$ with the following properties.  {\rm{(i)}} The elements of the Steiner quasigroup $Q$ are the points of
$\mathfrak{S}$. {\rm{(ii)}}
The restriction of the factor system of $\mathcal{Q}_f^{+}(\mathcal{Q}_f^{-})$ to
$(\mathfrak{S}\times \mathfrak{S})\setminus \{(x,x), x\in
\mathfrak{S}\}$ coincides with the orientation function of
$(\mathfrak{S},T)$. {\rm{(iii)}} $f(x,x)=1$, (respectively $f(x,x)=-1$) for
all $x\in Q$. \ed

\bt
Oriented Steiner quasigroups are flexible and have the cross inverse property.
\et

\bp
The statements follow from Remark 1 \cite{SS} p. 134. and Proposition 3.10 \cite{NS2} p. 767.
The bijection $\kappa$ that determines the left and the right cross inverse property has $\kappa:(a,\alpha)\mapsto (a,\alpha)$
for any $a\in Q$ and $\alpha\in Z_2$.
\ep

Since each element is its own cross inverse element, we obtain
\bc
Oriented Steiner quasigroups are semi-symmetric.
\ec

We note that the obtained classes of quasigroups differ from the
class of Mendelsohn quasigroups, since they are not idempotent. We have that $(a,-1)(a,-1)=(a,1)$
in $\mathcal{Q}_f^{+}$ and $(a,1)(a,1)=(a,-1)$ in $\mathcal{Q}_f^{-}$ for all $a\in Q$.

The oriented Steiner quasigroups do not
satisfy the left (right) alternative law (see Section \ref{sec4})
or the left (right) inverse property.
(The latter means that $\iota(a,\alpha)[(a,\alpha)(b,\beta)]=(b,\beta)$ with $\iota(a,\alpha)=(a,-\alpha)$ which
does not hold for $a=b$ and $\alpha = \beta$.) This shows a difference
with the case for oriented Steiner loops of exponent $4$. However,
the oriented Steiner quasigroups are flexible and have the cross inverse property,
like oriented Steiner loops
of exponent $2$ (see Theorem 5, \cite{SS}, p. 136.).

\subsection{Canonical oriented Steiner quasigroups}

We now turn back to the $f$-extensions of commutative quasigroups of order $3$. We saw in Section
\ref{ext3} that in order
to get a unique algebraic face of an oriented Steiner triple system as quasigroup $f$-extensions of
the Steiner quasigroup of order $3$ by Steiner quasigroups, we
can restrict our consideration to the case of the $f$-extensions of the group of order $3$ by Steiner quasigroups.

\bd A canonical oriented Steiner quasigroup $\mathcal{Q}_f$ is an
$f$-extension of the cyclic group $Z_3$ of order $3$
by the Steiner quasigroup $Q$ for
which there exists an oriented Steiner triple system
$(\mathfrak{S},T)$ with the following properties. {\rm{(i)}} The elements of the Steiner quasigroup $Q$
are the points of $\mathfrak{S}$. {\rm{(ii)}} The restriction of the factor system of $\mathcal{Q}_f$ to
$(\mathfrak{S}\times \mathfrak{S})\setminus \{(x,x), x\in
\mathfrak{S}\}$ coincides with the orientation function of
$(\mathfrak{S},T)$. {\rm{(iii)}} $f(x,x)=0$ for all $x\in Q$. \ed

\bt
Canonical oriented Steiner quasigroups satisfy the inverse property.
\et

\bp
$Q$ is commutative, has the inverse property and the bijection
$\iota:\mathcal{Q}_f \longrightarrow \mathcal{Q}_f:(a,\alpha)\mapsto (a,-\alpha)$
for all $a\in Q$ and $\alpha\in Z_3$ yields the right and the left inverse property of
$\mathcal{Q}_f$.
\ep

The canonical oriented Steiner quasigroups are not flexible.
This is deduced, with the help of Proposition 3.10 in \cite{NS2} p. 767, from the fact that
$$f(x,yx)f(y,x)\neq f(xy,x)f(x,y) \;\hbox{for all $x, y\in Q$.}$$

Next, the canonical oriented Steiner quasigroups are not idempotent.
This follows from the fact that
$$(a,\alpha)(a,\alpha)=(a,\alpha + \alpha)=(a,-\alpha)\;\hbox{ for all $a\in Q$, $\alpha\in Z_3$.}$$

Further, the canonical oriented Steiner quasigroups do not have the cross inverse property:
$$[(a,\alpha)(b,\beta)](x,\delta(a,\alpha))=(b,\beta)\;\hbox{ with $x=a$
and $\delta(a,\alpha)=(1-\alpha)$,}$$ and
$$
(y,\xi(a,\alpha))[(b,\beta)(a,\alpha)]=(b,\beta)\;\hbox{ with $y=a$
and $\xi(a,\alpha)=(-1-\alpha)$}.
$$
\noindent
In fact, both of them fail to hold in the case $a=b$, $\alpha=\beta$.

Finally, the canonical oriented Steiner quasigroups are not semi-symmetric:
$$(a,\alpha)[(b,\beta)(a,\alpha)]=(b,1+\alpha+\alpha+\beta)\neq (b,\beta).$$

\section{Near-associativity}

Canonical oriented and oriented Steiner quasigroups are not loops. Hence, to
measure their near-associativity we can consider Belousov's generalization of the notion of nuclei (the orbits of the
groups of regular permutations of the quasigroup). Cf. \cite{B} p. 22--25.

A bijection $\lambda:Q \to Q$ is a {\it left-regular permutation}
or a {\it right-regular permutation} of $(Q,\cdot)$, if for all $x,y\in Q$
one has $\lambda (xy)=\lambda (x)\cdot y$ or $\rho (xy)=x\cdot\rho (y)$,
respectively.
If $\lambda$ is a left-regular permutation then
$\lambda = \lambda_{\lambda(x)}\lambda^{-1} _x$ for all $x\in Q$.
Similarly, if $\rho$ is a right-regular permutation then $\rho =
\rho_{\rho(x)}\rho^{-1} _x$ for all $x\in Q$.
Hence the left-regular permutations form a subgroup $\Lambda(Q)$ of the left multiplication group
$G_l(Q)$, and the right-regular permutations form a subgroup $\text{R}(Q)$ of
the right multiplication group $G_r(Q)$. If $\Lambda(Q)$ or $\text{R}(Q)$ consists only of the identity
map of $Q$, then we say that the group of left-regular permutations or the group of right-regular
permutations is {\it trivial}.

Steiner quasigroups are idempotent; consequently, they have only tri\-vi\-al right regular permutations.
According to Remark 3.1 in \cite{S}, p. 113., the orbits of the group of right regular permutations
of Steiner quasigroups are contained in the congruence classes of the $f$-extensions. Since
we extend by idempotent quasigroups, these equivalence classes are normal sub-quasigroups
of (canonical) oriented Steiner quasigroups.

\bs
The group of the right regular permutations and the group of the left regular permutations, of
an oriented Steiner quasigroup, are both isomorphic to the group $Z_2$. The
group of right regular permutations and the group of left regular permutations, of
a canonical oriented Steiner quasigroup, are both isomorphic to the group $Z_3$.
\es

\bp
The claims are deduced from Theorem 3.1 and Theorem 4.2 in \cite{S}, p. 113. and p. 115.
\ep

\br
In the spirit of this generalization of the nucleus,
(canonical) oriented Steiner quasigroups are
right nuclear $f$-extensions. This is true because the orbits of the group of right regular
permutations of the
oriented Steiner quasigroup $R(\mathcal{Q}_f^{+})$ (\,$R(\mathcal{Q}_f^{-})$) and
of the canonical oriented Steiner quasigroup $R(\mathcal{Q}_f)$ coincide with the
congruence classes of the $f$-extensions $\{(a,\alpha): \alpha \in Z_2\}$ and
$\{(a,\alpha): \alpha \in Z_3\}$, respectively.

Also, (canonical) oriented Steiner quasigroups are left nuclear extensions, because
the same argument is valid for the
orbits of the group of left regular permutations.

Thus, (canonical) oriented Steiner quasigroups are nuclear $f$-extensions.
\er

\section{Oriented Steiner triple systems in cryptology}

It is commonly recognized that cryptology consists of two parts: {\it cryptography} - used for ``defence", i.e., for constructing ciphers, - and {\it cryptanalysis} - used for ``attacks", i.e., for developing methods on breaking ciphers.
(The latter is often treated as a kind of an ``art" based on case-by-case actions.) A cipher is a device, or a method, of information encryption, used with a purpose of enhancing security. A {\it private key} is a ``hidden" collection of parameters
of a cipher made known to the recipient but not divulged to the public. A number of constructions of error detecting and error correcting codes, cryptographic algorithms and enciphering systems have been constructed with the help of some associative algebraic structures. However, there exists also a possibility of using non-associative structures such as quasigroups and neo-fields. Sometimes codes and ciphers based on non-associative structures display better qualities than codes and ciphers based on associative structures. On the other hand, efficiency of applications of quasigroups in both cryptography and cryptanalysis is based on the fact that quasigroups can be interpreted as generalized ``permutations". (The number of quasigroups of order $n$ exceeds $n!\cdot(n-1)!\cdot \; ... \; \cdot 2!\cdot1!$.) The ratio of the number of groups to the number of quasigroups of a given size tends to zero as the size tends to infinity.
A survey of results about quasigroups in cryptology can be found in \cite{Sh} and in the references therein.

\bigskip
Now, we propose an encryption algorithm based on quasigroups using their extensions.
Let $K$, $Q$ be arbitrary quasigroups, and consider the Schreier-type quasigroup extension
$\mathcal{Q}$ defined on $Q\times K$ by the operation:
$$
(a,\alpha)\circ (b,\beta) = (ab, f(a,b) \alpha^{G(b)}\beta),
$$
where $f:Q\times Q\longrightarrow K$, and $G:Q \longrightarrow Aut(K)$.

The message is a string $q=(q_1, q_2, ... ,q_n)$ consisting of elements of $Q$. The public key
is formed by strings $k=(k_1, k_2 ..., k_n)$, $(k_i\in K)$ and $c=(c_1, ... ,c_n)$, $c_i\in \mathcal{Q}$. The private key consists of function $f$ and map $G$. The enciphered message is a string
$a=(a_1=(q_1,k_1)\circ c_1, a_2 =(q_2,k_2)\circ c_2, ... ,a_n=(q_n,k_n)\circ c_n)$ consisting of elements of $\mathcal{Q}$.
The sender broadcasts both the ciphertext $a$ and the keys $k$, and $c$.

Quasigroup $\mathcal{Q}$ (i.e., its multiplication table) is also distributed as a public part of the key. The security lies in the function $f$ and in the map $G$. For sufficiently large quasigroups
the cryptanalyst will face a difficult task to decode the ciphertext without knowing $f$ and $G$.
In fact, a possible interceptor may have $a$, $k$ and $c$, but must
find out the original message $q$; and even possible successes with deciphering
past messages may not help in future in the case of large quasigroups $K,Q$.

The intended receiver of course knows $f$, and $G$, gets $a$, $k$, $c$.

As an example, we apply a simplified version of this algorithm, to show how oriented Steiner
triple systems can occur in enciphering, using canonical oriented Steiner quasigroups.

As before, let $(\mathfrak{S}, T)$ be an oriented Steiner triple system, $Q$ be a Steiner quasigroup
formed by points of $(\mathfrak{S}, T)$. Then $\mathcal{Q}_f$ is the associated canonical oriented Steiner quasigroup (i.e., the $f$-extension of $Z_3$ by $Q$ with the factor system $f$, determined by the orientation function of the oriented Steiner triple system $(\mathfrak{S}, T)$).

A message $q$ is given as a sequence $q_1, ... ,q_n$ of points of the oriented Steiner triple system, (elements of the quasigroup $Q$), and public keys are sequences $k_1, ... ,k_n\in Z_3$ and $c_1, ... ,c_n\in \mathcal{Q}_f$. The ciphertext is $a_1 = (q_1,k_1)\circ c_1, ... ,a_n = (q_n,k_n)\circ c_n$.

The sender transmits the enciphered message and keys. Orientations of blocks are known only by the
intended receiver.

The effectiveness of this algorithm is related to the fact that the blocks of the
Steiner triple system are oriented independently. Thus, from a Steiner triple systems with $n$ points and $b=\frac{n(n-1)}{6}$ blocks one can construct $2^b$ oriented Steiner triple systems.

The same algorithm is working in the case of oriented Steiner quasigroups.

\vskip 10pt
{\bf Acknowledgement}
\vskip 10pt
\noindent
The author has been supported by FAPESP Grant - process No 11/51845-5,
and expresses her gratitude to IMS, University of S\~{a}o Paulo,
Brazil, for the warm hospitality.

\vskip 10pt
Izabella Stuhl\\
Institute of Mathematics and Statistics\\
University of S\~{a}o Paulo\\
05508-090 S\~{a}o Paulo, SP, Brazil \vspace{1mm}\\
University of Debrecen\\
H-4010 Debrecen, Hungary\\
{\it E-mail}: {\it {}izabella@ime.usp.br}\\

\end{document}